\documentclass{article}

\pdfoutput=1

\usepackage{amsmath}
\usepackage{amsthm}
\usepackage{amsfonts}
\usepackage{amssymb}
\usepackage{color}
\usepackage[colorlinks=true, pdfstartview=FitV, linkcolor=blue, citecolor=blue, urlcolor=blue]{hyperref}
\newtheorem{thm}{Theorem}[section]
\newtheorem{cor}[thm]{Corollary}
\newtheorem{lem}[thm]{Lemma}

\theoremstyle{definition}

\title          {Bounds on the roots of the Steiner Polynomial}
\author         {Madeleine Jetter}

\begin {document}

\maketitle

\begin{abstract} We consider the Steiner polynomial of a $C^2$ convex body $K \subset \mathbb{R}^n$. Denote by $\rho_{\min}$ the minimum value of the principal radii of curvature of $\partial K$ and by $\rho_{\max}$ their maximum. When $n \leq 5$, the 
real parts of the roots are bounded above by $-\rho_{\min}$ and below by $-\rho_{\max}$. These bounds are valid for any $n$ such that all of the roots of the Steiner polynomial of every convex body in $\mathbb{R}^n$ lie in the left half-plane. 
\end{abstract}

\section{Introduction}

Let $K \subset \mathbb{R}^n$ be a (compact) convex body and let $B$ denote the unit ball in $\mathbb{R}^n$. We form the outer parallel body $K+tB$ by taking the Minkowski sum of $K$ and a ball of radius $t > 0$, that is:  
$$K+tB = \left\{ \vec{v} + t \vec{u} \mid \vec{v} \in K, \vec{u} \in B \right\}.$$
Thinking of the outer parallel body as the result of the unit-speed outward normal flow applied to $K$ at time $t$ makes it relevant to applied problems such as combustion \cite{GrOs}. The volume of $K + tB$ can be written as a polynomial of degree $n$, the Steiner polynomial \cite{Sch}:
$$V_{K + tB} = \sum_{i = 0}^n {n \choose i} V(K^{n-i}, B^i) t^i$$
where the coefficient $V(K^{n-i}, B^i)$ is the mixed volume of $n-i$ copies of $K$ and $i$ copies of the unit ball. We will adopt the notation $S_K(t) = V_{K+tB}$ for the Steiner polynomial of $K$ in the variable $t$.

In two dimensions, consideration of the roots of the Steiner polynomial leads to a Bonnesen-style inequality. When $K$ is a convex planar region, with area $A_K$ and perimeter $L_K$, we have 
$$S_K(t) = A_K + L_K t + \pi t^2.$$
Since the discriminant of the Steiner polynomial in two dimensions is $L_K^2 - 4 \pi A_K$,
we see that the isoperimetric inequality for $K$ is equivalent to the fact that $S_K (t)= 0$ has (one double or two single) real roots. Moreover, since $S_K$ gives the area of the region $K + tB$, the roots must also be negative when $A_K > 0$. 

Furthermore, it is known that
\begin{thm}Let $K$ be a strictly convex region which is not a disc. Let $R_i = \sup \{ r \mid  \textrm{a translate of} \ r B \subset K \}$ be the inradius of $K$, and let $R_e = \inf \{ r \mid  \textrm{a translate of} \ K \subset r B \}$ be the outradius. Let $\rho_{\min}$ and $\rho_{\max}$ denote the minimum and maximum values of the radius of curvature of $K$. If the roots of $S_K$ are $t_1 < t_2$, then
\begin{equation}\label{PlanarComparison}-\rho_{\max} < t_1 < -R_e < -\frac{L_K}{2\pi} < -R_i < t_2 < -\rho_{\min}.\end{equation}
\end{thm}

When $K$ is a disc, then all of the above quantities are equal, giving a version of Bonnesen's inequality. Green and Osher provide a proof in \cite{GrOs}.

Teissier \cite{Te:B}, working in the setting of ample divisors on algebraic varieties, posed the following problems aimed at generalizing the appealing state of affairs in the planar case. Suppose a convex body  $K \subset \mathbb{R}^n$ is given and that the roots of $S_K$ have real parts $r_1 \leq r_2 \leq \dots \leq r_n$. 

\begin{itemize}
\item[P1.] Is $S_K$ stable (i.e. do all the roots lie in the left half-plane)? 
\item[P2.] Let $R_i$ indicate the inradius of $K$, that is, the largest real number $s$ such that a translate of $sB$ is contained in $K$. Does the inequality $-R_i \leq r_n$ hold?
\end{itemize}

By the Routh-Hurwitz stability criterion and the Aleksandrov-Fenchel inequalities \cite{Sch}, we know that $S_K$ is stable for all convex bodies  $K \subset \mathbb{R}^n$ provided that $n \leq 5$. On the other hand, Cifre and Henk construct an example in \cite{CiHe} to show that $S_K$ need not be stable when $K \subset \mathbb{R}^{15}$. Less is known about the inradius bound. However, in those cases where Teissier's first problem has an affirmative answer, we can prove a generalization of the extreme upper and lower bounds in inequality (\ref{PlanarComparison}) relatively easily.

\begin{thm}\label{mainthm}Assume that in $\mathbb{R}^n$, $S_K$ is stable for every convex $K \subset \mathbb{R}^n$. Let $K \subset \mathbb{R}^n$ be a $C^2$ convex body, and suppose that the roots of $S_K$ have real parts $r_1\leq r_2 \leq \dots \leq r_n$. Denote by $\rho_{\min}$ and $\rho_{\max}$ the minimum and maximum values of the principal radii of curvature of $K$. Then 

 $$-\rho_{\max} \leq r_1 \leq \dots \leq r_n \leq -\rho_{\mathrm{min}}. $$
\end{thm}

\section{Technical Background}\label{TechPrereq}
\subsection{The Steiner Polynomial}
A general reference for this section is Schneider's volume \cite{Sch}.
The fundamental tool for what follows is the support function $p_K: \mathbb{R}^n \rightarrow \mathbb{R}$ of a convex body $K \subset \mathbb{R}^n$, defined as follows:

$$p_K(\vec{x}) = \sup \{ \vec{x} \cdot \vec{v} \mid \vec{v} \in K \},$$

where $\cdot$ denotes the standard inner product. Because of the homogeneity of the support function, $p_K$ is determined by its restriction to the unit sphere. Thus we frequently treat $p_K$ as a function on $S^{n-1}$.

A particularly important feature of $p_K$ is the way in which it carries information about the curvature of $\partial K$ when the boundary satisfies certain smoothness conditions. When $K$ (and thus $p_K$) is $C^2$, we consider the Hessian matrix $H(p_K)$. Given $\omega \in S^{n-1}$, we choose a basis $\left\{ e_1, \ldots , e_n \right\}$ where $\left\{ e_1, \ldots , e_{n-1} \right\}$ is an orthonormal basis for $TS^{n-1}_{\omega}$ and $e_n = \omega$. One can show using homogeneity \cite{Sch} that the eigenvalues of $H(p_K(\omega))$ computed with respect to this basis are $0$ and the principal radii of curvature of $K$ at $\omega$, which we denote $\rho_1, \dots , \rho_{n-1}$. The restriction of the Hessian to $TS^{n-1}_{\omega}$, which we write as $\overline{H}(p_K(\omega))$, has eigenvalues $\rho_1, \ldots , \rho_{n-1}$. 

Since the infinitesimal element of area on $\partial K$ is the product of the principal radii of curvature, 
we may write the volume of $K$ equivalently as
 $$V_K = \frac{1}{n} \int_{S^{n-1}} p_K \rho_1 \cdot \dots \cdot \rho_{n-1}  \ d\omega$$
or
$$V_K = \frac{1}{n} \int_{S^{n-1}} p_K \det \overline{H}(p_K) \ d\omega$$
where $\overline{H}$ denotes the Hessian matrix computed with respect to an orthonormal frame for $TS^{n-1}$.

Applying this formula to $K + t B$, noting that $p_{K + tB} = p_K + t$, we have
\begin{equation}\label{SteinerPolyIntegral}
S_K = V_{K + tB} = \frac{1}{n} \int_{S^{n-1}} (p_K + t) \det \left( \overline{H}(p_K) + t I \right) \ d\omega .
\end{equation}

The integrand above is a polynomial of degree $n$ in $t$. We can isolate the coefficient of each $t^i$ using the Minkowski integral formulas (\cite{Sch}, p. 291) to obtain an integral expression for $V(K^{n-i}, B^i)$. 
\begin{align*}
V(K^{n-i}, B^i) &= \frac{1}{n} \int_{S^{n-1}} s_{n-i}(\rho_1, \dots , \rho_{n-1}) \ d \omega \\
 &=  \frac{1}{n} \int_{S^{n-1}} p_K s_{n-i-1}(\rho_1, \dots , \rho_{n-1} ) \ d \omega ,
 \end{align*}
where $s_j$ is the normalized $j^{\mathrm{th}}$ elementary symmetric function in $\rho_1, \dots , \rho_{n-1}$ (ie. ${n-1 \choose j} s_j$ is the usual  $j^{\mathrm{th}}$ elementary symmetric function).

\subsection{Minkowski Subtraction}
The proof of theorem \ref{mainthm} will also rely on the concept of Minkowski subtraction. Given convex bodies $K, L \subset \mathbb{R}^n$, the \emph{Minkowski difference} of $K$ and $L$ is
$$K \sim L = \left\{\vec{v} \in \mathbb{R}^n \mid L + \vec{v} \subset K \right\}. $$ 

We may think of $K \sim L$ as the intersection of all translates of $K$ by opposites of vectors in $L$. If $K$ and $L$ are both convex, then $K \sim L$ is as well, but the operations of Minkowski sum and difference are not  inverse to one another. Although $(K + L) \sim L = K$ holds for any convex bodies $K$ and $L$, $(K \sim L) + L = K$ only when there exists a convex body $M$ such that $L + M = K$. In this case we say that $L$ is a Minkowski summand of $K$, and $(K \sim L) + L = M + L = K.$ 

Specializing to a situation relevant to the proof, when we know that $cB$ is a Minkowski summand of $K$, we have that $(K \sim cB) + cB = K$ and it follows that $p_{K \sim cB} = p_K - c$. This allows us to compute $S_{K \sim cB}$ fairly easily using Equation (\ref{SteinerPolyIntegral}). 

We will make use of the following lemma appearing in \cite{Ma} and \cite{Sch}, which gives a condition under which $L$ is a Minkowski summand of $K$.

\begin{lem}\label{FreelySliding} Suppose $K, L \subset \mathbb{R}^n$ are convex. If the maximum of all the principal radii of curvature of $L$ is bounded above by the minimum of the principal radii of curvature of $K$ at each $\omega \in S^{n-1}$, then $L$ is a Minkowski summand of $K$ -- i.e. there is a convex body $M$ such that $L + M = K$.\end{lem}

\section{Proof of Theorem \ref{mainthm}}
We first establish the upper bound, which is the easier of the two. Since $K$ is convex, each $\rho_i \geq 0$. We may assume that $K$ is $C^2_+$, (in other words the principal radii of curvature are all strictly positive and hence $\rho_{\min} > 0$) since otherwise there is nothing to prove. 
If $0 \leq c \leq \rho_{\mathrm{min}}$, then let $K' = K \sim cB$. $cB$ is a Minkowksi summand of $K$ by Lemma \ref{FreelySliding}, so 

$$S_{K'}(t) = \frac{1}{n} \int_{S^n-1} (p_K - c + t)  \mathrm{det} \left( \overline{H}(p_K) + (-c + t)I \right) d\omega= S_K(t - c).$$ 

The roots of $S_{K'}$ have real parts $r_i + c$, so the stability assumption implies  $r_i + c < 0$, hence $r_i < -c$ for any $c \leq \rho_{\mathrm{min}}$. Letting $c = \rho_{\mathrm{min}}$ yields the claimed upper bound.

Turning to the lower bound, let $c \geq \rho_{\mathrm{max}}$. Then $K$ is a Minkowski summand of $cB$ and we write $K' = cB \sim K$. Writing $p_K$ for the support function of $K$, we have $p_{K'} = c - p_K$. Expanding the Steiner polynomial of $K'$,
$$\frac{1}{n} \int_{S^n-1} (-p_K + c + t) \mathrm{det} \left( \overline{H}(-p_K) + (c + t)I \right) d\omega,$$
in the case $n=3$ we have $S_{K'} = - \left( V_K - A_K (c+t) + H_K (c+t)^2 - V_B (c+t)^3 \right)$, and in general $S_{K'} = (-1)^n S_K(-t-c)$. The roots of $S_{K'}$ have real parts $-(r + c)$, so by stability $-r - c < 0$ and we conclude that $-c < r$. The lower bound follows by taking $c =  \rho_{\mathrm{max}}$.

\begin{cor}The real parts of the roots of $S_K$ are bounded by $-\rho_{\min}$ and $-\rho_{\max}$ for any $C^2$ convex body $K \subset \mathbb{R}^n$ where $n \leq 5$.\end{cor}

\begin{proof} It is known \cite{Te:B} that for $n \leq 5$, $S_K$ is stable for every convex body $K \subset \mathbb{R}^n$. This follows from the Routh-Hurwitz stability criterion and the Aleksandrov-Fenchel inequalities.
\end{proof}

\bibliography{MJbib}
\bibliographystyle{hplain}

\end {document}